\documentclass{article}
\usepackage{amssymb}
\usepackage{amsfonts}
\usepackage{graphicx}
\usepackage{amsmath}

\newtheorem{theorem}{\sc{Theorem}}[section]

\newtheorem{definition}[theorem]{\sc{Definition}}
\newtheorem{example}[theorem]{\sc{Example}}
\newtheorem{lemma}[theorem]{\sc{Lemma}}

\newtheorem{remark}[theorem]{\sc{Remark}}
\newenvironment{proof}[1][Proof]{\textbf{#1.} }{\ \rule{0.5em}{0.5em}}
\input{tcilatex}

\begin{document}

\title{Generic Lie Color Algebras}
\footnotetext{Mathematics Subject Classification: Primary 17B75, Secondary 17B35}
\author{Kenneth L. Price \\
Department of Mathematics\\
University of Wisconsin\ Oshkosh\\
Oshkosh, WI 54901-8631\\
E-mail: pricek@uwosh.edu}
\maketitle

\begin{abstract}
We describe a type of Lie color algebra, which we call generic, whose
universal enveloping algebra is a domain with finite global dimension.
Moreover, it is an iterated Ore extension. We provide an application and
show Gr\"{o}bner basis methods can be used to study universal enveloping
algebras of factors of generic Lie color algebras.
\end{abstract}

\section{Introduction}

Throughout $k$ denotes a field of characteristic different from 2. Iterated
Ore extensions play an important role in this work. All of the Ore
extensions we introduce are of the form $A\left[ \theta ;\sigma ,\delta %
\right] $ where $A$ is a $k$-algebra, $\sigma $ is an automorphism, and $%
\delta $ is a $\sigma $-derivation.\ In this case if $A$ is Noetherian, a
domain, or has finite right global dimension, then $A\left[ \theta ;\sigma
,\delta \right] $ has that corresponding property as well. These facts are
well known (see \cite[Theorem I.2.9 and Theorem VII.5.3 ]{McRob}) and will
be invoked as needed without further comment.

We follow the notation on Lie color algebras in \cite[Sections 1, 2, and 6]%
{Prime}.\ Let $k^{\times }=k\backslash \left\{ 0\right\} $ denote the group
of units of $k$.

\begin{definition}
Let $G$ be an abelian group. A map $\varepsilon :G\times G\rightarrow
k^{\times }$ is called a \emph{skew-symmetric bicharacter} on $G$ if it
satisfies (1) and (2) below, for any $f,g,h\in G$.

\begin{enumerate}
\item $\varepsilon (f,g+h)=\varepsilon (f,g)\varepsilon (f,h)$ and $%
\varepsilon (g+h,f)=\varepsilon (g,f)\varepsilon (h,f)$

\item $\varepsilon (g,h)\varepsilon (h,g)=1$
\end{enumerate}
\end{definition}

All gradings are with respect to an additively written abelian group. We
denote the degree of a nonzero homogeneous element $x$ by $\partial x$. We
want to avoid statements like \textquotedblleft $\partial x=g$ or $x=0$%
\textquotedblright\ so we often write $\partial x=g$ to handle both cases.
If $x$ and $y$ are homogeneous and $\varepsilon $ is a skew-symmetric
bicharacter, then we shorten our notation by writing $\varepsilon \left(
x,y\right) $ instead of $\varepsilon \left( \partial x,\partial y\right) $.

\begin{definition}
A $(G,\varepsilon )$\emph{-Lie color algebra} is a $G-$graded vector space $%
\mathcal{L}$ equipped with a graded bilinear map $\langle ,\rangle :\mathcal{%
L}\times \mathcal{L}\rightarrow \mathcal{L}$, called the \emph{bracket} of $%
\mathcal{L}$, which satisfies the following for any homogeneous $x,y,z\in 
\mathcal{L}$.

\noindent $%
\begin{array}[t]{cl}
\langle x,y\rangle =-\varepsilon (x,y)\langle y,x\rangle & \varepsilon \text{%
-skew-symmetry} \\ 
\varepsilon (z,x)\langle x,\langle y,z\rangle \rangle +\varepsilon
(y,z)\langle z,\langle x,y\rangle \rangle +\varepsilon (x,y)\langle
y,\langle z,x\rangle \rangle =0 & \varepsilon \text{-Jacobi identity}%
\end{array}%
$
\end{definition}

A skew-symmetric bicharacter $\varepsilon $ satisfies $\varepsilon (g,g)=\pm
1$ for each $g\in G$. Note that $G_{+}=\{g\in G:\varepsilon (g,g)=1\}$ is a
subgroup such that $[G:G_{+}]\leq 2$. Set $G_{-}=\{g\in G:\varepsilon
(g,g)=-1\}=G\backslash G_{+}$. For any $G$-graded vector space $V$ we set $%
V_{\pm }=\oplus _{g\in G_{\pm }}V_{g}$.

We may view Lie superalgebras and graded Lie algebras, which are graded over 
$\mathbb{Z}_{2}$ and $\mathbb{Z}$, respectively, as special types of Lie
color algebras.\ The `absolutely torsion free' condition was introduced by
R. B\o gvad to find graded Lie algebras with finite global dimension (see 
\cite[Theorem 1]{Bogvad}). It is well-known that the universal enveloping
algebra of a finite dimensional Lie superalgebra may have infinite global
dimension (see \cite[Proposition 5]{Behr}). M. Aubry and J.-M. Lemaire have
shown that the universal enveloping algebra of an absolutely torsion free
graded Lie algebra (or Lie superalgebra) is a domain (see \cite{AL}).

Now consider a finite dimensional Lie color algebra $\left( \mathcal{L}%
,\left\langle ,\right\rangle \right) $. Its universal enveloping algebra $%
U\left( \mathcal{L}\right) $ need not even be semiprime (see \cite[Example
2.9]{Prime}). Conditions for $U\left( \mathcal{L}\right) $ to be semiprime
or prime are provided in \cite[Thoerem 2.5]{Prime}, which extends the
analogous result, \cite[Theorem 1.5]{Bell}, for Lie superalgebras. When $%
U\left( \mathcal{L}\right) $ has finite global dimension, then it must equal 
$\dim \mathcal{L}_{+}$ by \cite[Theorem 3.1]{Homolog}, which extends the
analogous result, \cite[Proposition 2.3]{KK}, for Lie superalgebras.

We do not know how to extend Aubry and Lemaire's or B\o gvad's theorems to
Lie color algebras (see \cite{AL}\ and \cite{Bogvad}). However we believe it
may be achieved through the use of generic Lie color algebras, which are
defined in section \ref{U(X) section}.

We handle the case with $\dim \mathcal{L}_{-}\leq 2$ in section \ref{Global
Dim Sec}. In this case Theorem \ref{Global Dim Thm} provides conditions for $%
U\left( \mathcal{L}\right) $ to be a domain with finite global dimension.
Example \ref{Color Example} shows this is possible even when there exists $%
x\in \mathcal{L}_{-}$ such that $\left\langle x,x\right\rangle =0$. Results
of Behr, B\o gvad, and Aubry and Lemaire (discussed above) imply that this
is not possible for Lie superalgebras or graded Lie algebras. Our proof uses
generic Lie color algebras.

In section \ref{Grob Method} we explain how to find the Gr\"{o}bner basis of
an ideal generated by positive elements of a generic Lie color algebra. We
refer the reader to \cite{Effective} for background on Gr\"{o}bner bases, to 
\cite{Primality Test} for a Gr\"{o}bner basis test to determine if an ideal
of an iterated Ore extension is completely prime, and to \cite{Homological}
for a Gr\"{o}bner basis method to calculate projective dimension.

\section{Generic Lie Color Algebras \label{U(X) section}}

For a Lie color algebra $\left( \mathcal{L},\left\langle ,\right\rangle
\right) $ and a linear subspace $V$\ of $\mathcal{L}_{-}$ we let $%
\left\langle V,V\right\rangle $ denote the linear subspace of $\mathcal{L}%
_{+}$ generated by brackets between elements of $V$. If $n=\dim V$ then it
is easy to see that $\dim \left\langle V,V \right\rangle \leq \frac{1}{2}%
n\left( n+1\right) $.

\begin{definition}
\label{generic Lie color algebra}A Lie color algebra $\left( \mathcal{X}%
,\left\langle ,\right\rangle \right) $ is called \emph{generic} if $\dim 
\mathcal{X}_{+}=\frac{1}{2}m\left( m+1\right) $, where $m=\dim \mathcal{X}%
_{-}<\infty $, and $\mathcal{X}_{+}=\left\langle \mathcal{X}_{-},\mathcal{X}%
_{-}\right\rangle $ is color central, that is, $\left\langle
x,y\right\rangle =0$ for all $x\in \mathcal{X}_{+}$ and $y\in \mathcal{X}$.
\end{definition}

\begin{remark}
\label{generic subalgebra}Let $V$ be a graded subspace of $\mathcal{X}_{-}$.
We set $\mathcal{L}_{-}=V$ and $\mathcal{L}_{+}=\left\langle
V,V\right\rangle $. Then it is easy to show the sub Lie color algebra $%
\mathcal{L}$ of $\mathcal{X}$ is also a generic Lie color algebra.
\end{remark}

\begin{lemma}
\label{Images of Generic}Let $\left( \mathcal{L},\left\langle ,\right\rangle
\right) $ be a finite dimensional Lie color algebra such that $\mathcal{L}%
_{+}=\left\langle \mathcal{L}_{-},\mathcal{L}_{-}\right\rangle $ is color
central. Then there is a generic Lie color algebra $\mathcal{X}$ with $\dim 
\mathcal{X}_{-}=\dim \mathcal{L}_{-}$ and a surjective Lie color algebra
homomorphism $\psi :\mathcal{X}\rightarrow \mathcal{L}$.
\end{lemma}

The proof of Lemma \ref{Images of Generic} is straightforward.

Given a $G$-graded algebra $A$ we write $A^{\varepsilon }\left[
t_{1},t_{2},\ldots ,t_{n}\right] $ to denote the color polynomial algebra
over $A$ in $n$ homogeneous variables. It is $G$-graded and isomorphic to an
iterated Ore extension of $A$ (see \cite{BMPZ} for details).

\begin{theorem}
\label{Generic Ore}Suppose $\mathcal{X}$ is generic and $x_{1},\ldots ,x_{m}$
form a homogeneous basis of $\mathcal{X}_{-}$. Set $g_{i}=\partial x_{i}$
and $t_{i,j}=\frac{1}{2}\left\langle x_{i},x_{j}\right\rangle $ for $1\leq
i<j\leq n$. Then $U\left( \mathcal{X}\right) $ is isomorphic to the iterated
Ore extension 
\begin{equation}
U\left( \mathcal{X}\right) \cong S\left[ x_{1};\alpha _{1}\right] \left[
x_{2};\alpha _{2},\delta _{2}\right] \ldots \left[ x_{m};\alpha _{m},\delta
_{m}\right]  \label{IteratedOre}
\end{equation}%
with $S=k^{\varepsilon }\left[ t_{i,j}:1\leq i<j\leq m\right] $ such that
for all $a\leq m$ and $1\leq b<c\leq m$

\begin{enumerate}
\item $\alpha _{a}\left( t_{b,c}\right) =\varepsilon \left(
g_{a},g_{b}+g_{c}\right) t_{b,c}$,

\item $\delta _{a}\left( t_{b,c}\right) =0$,

\item $\alpha _{c}\left( x_{b}\right) =\varepsilon \left( g_{c},g_{b}\right)
x_{b}$,

\item and $\delta _{c}\left( x_{b}\right) =-2\varepsilon \left(
g_{c},g_{b}\right) t_{b,c}=\left\langle x_{c},x_{b}\right\rangle $.
\end{enumerate}
\end{theorem}

Theorem \ref{Generic Ore} can be proved in the same way as \cite[Theorem
II.3.1]{Le Bruyn}. \ In fact the definition of generic Lie color algebra was
motivated by the treatment of generic Clifford algebras in \cite[Chapter 2]%
{Le Bruyn}.

\section{An Application \label{Global Dim Sec}}

We assume $k\ $is algebraically closed throughout this section.

\begin{lemma}
Let $\left( \mathcal{L},\left\langle ,\right\rangle \right) $ be a Lie color
algebra such that $\mathcal{L}_{+}$\ is color central and $\dim \mathcal{L}%
_{-}=\dim \mathcal{L}_{+}=2$. Suppose there does not exist homogeneous $x\in 
\mathcal{L}_{-}$ such that $\left\langle x,x\right\rangle =0$. Then $U\left( 
\mathcal{L}\right) \cong k\left[ \theta _{1}\right] \left[ \theta
_{2};\sigma _{q}\right] $ where $q$ is a nonzero scalar, $\sigma _{q}\left(
\theta _{1}\right) =q\theta _{1}$, and either $\theta _{1},\theta _{2}$ are\
homogeneous or $\theta _{1}+\theta _{2},\theta _{1}-\theta _{2}$ are
homogeneous and $q=-1$.\label{Surprise}
\end{lemma}

\proof%
In view of Lemma \ref{Images of Generic} we may pass to the case that $%
\mathcal{L}\cong \mathcal{X}/K$ where $K=\ker \psi $ is a homogeneous linear
subspace of $\mathcal{X}_{+}$ with $\dim K=1$.

\begin{description}
\item[Step 1.] Choose nonzero homogeneous $t\in K$. There is a homogeneous
basis $\left\{ x_{1},x_{2}\right\} $ of $\mathcal{X}_{-}$ such that either $%
t=\lambda _{1}\left\langle x_{1},x_{1}\right\rangle +\lambda
_{2}\left\langle x_{2},x_{2}\right\rangle $ or $t=\left\langle
x_{1},x_{2}\right\rangle $. Set $g_{1}=\partial x_{1}$ and $g_{2}=\partial
x_{2}$.
\end{description}

Let $y_{1},y_{2}$ form a homogeneous basis of of $\mathcal{X}_{-}$. Then we
may write $t$ as in equation \ref{t form} for some $\mu _{1},\mu _{2},\mu
_{3}\in k$.%
\begin{equation}
t=\mu _{1}\left\langle y_{1},y_{1}\right\rangle +\mu _{2}\left\langle
y_{2},y_{2}\right\rangle +\mu _{3}\left\langle y_{1},y_{2}\right\rangle
\label{t form}
\end{equation}%
If $\mu _{3}=0$ or $\mu _{1}=\mu _{2}=0$ then step 1 follows immediately.\
We pass to the case $\mu _{2}\neq 0$ and $\mu _{3}\neq 0$ by relabeling $%
y_{1}$ and $y_{2}$ if necessary. This implies $\partial t=2\partial
y_{2}=\partial y_{1}+\partial y_{2}$, which yields $\partial y_{1}=\partial
y_{2}$. Set $\lambda _{1}=\mu _{1}-\left( 2^{-1}\mu _{3}\right) ^{2}\left(
\mu _{2}\right) ^{-1}$, $x_{1}=y_{1}$, $\lambda _{2}=\mu _{2}$, and $%
x_{2}=y_{2}+\mu _{3}\left( 2\mu _{2}\right) ^{-1}y_{1}$. A straightforward
calculation shows $t=\lambda _{1}\left\langle x_{1},x_{1}\right\rangle
+\lambda _{2}\left\langle x_{2},x_{2}\right\rangle $.

\begin{description}
\item[Step 2.] Set $T=U\left( \mathcal{X}\right) $. Then $T\cong k\left[
t_{12}\right] \left[ x_{1};\sigma _{1}\right] \left[ x_{2};\sigma
_{2},\delta _{2}\right] $, defined as in Theorem \ref{Generic Ore}, and $%
U\left( \mathcal{L}\right) \cong T/\left( t\right) $.
\end{description}

The last statement follows from universal properties of enveloping algebras
and Ore extensions.

\begin{description}
\item[Step 3.] If $t=\left\langle x_{1},x_{2}\right\rangle $ the lemma holds
with $q=\varepsilon \left( g_{1},g_{2}\right) $. In this case $\theta _{1}$
and $\theta _{2}$ are homogeneous.
\end{description}

This is easy since $t=2t_{12}$ so $U\left( \mathcal{L}\right) \cong T/\left(
t_{12}\right) $.

\begin{description}
\item[Step 4.] Suppose $t=\lambda _{1}\left\langle x_{1},x_{1}\right\rangle
+\lambda _{2}\left\langle x_{2},x_{2}\right\rangle $. The lemma holds with $%
q=-1$. If $\partial x_{1}=\partial x_{2}$ then $\theta _{1}$ and $\theta
_{2} $ are homogeneous. Otherwise $\theta _{1}+\theta _{2}$ and $\theta
_{1}-\theta _{2}$ are homogeneous.
\end{description}

If $\lambda _{i}=0$ for $i=1$ or $i=2$ then $\left\langle x,x\right\rangle
=0 $ with homogeneous $x=\psi \left( x_{i}\right) \in \mathcal{L}_{-}$. We
assumed this could not happen. Thus $\lambda _{1}\neq 0$ and $\lambda
_{2}\neq 0$ which implies $\partial t=2g_{1}=2g_{2}$ and $\varepsilon \left(
g_{1},g_{2}\right) ^{2}=1$. By replacing $x_{1}$ and $x_{2}$ by appropriate
scalar multiples, if necessary, we may assume $\lambda _{1}=1$ and $\lambda
_{2}=-1$.

Set $S=k\left[ u,z_{1}\right] \left[ z_{2};\sigma ,\delta \right] $, where $%
\sigma \left( z_{1}\right) =-z_{1}$, $\sigma \left( u\right) =u$, $\delta
\left( z_{1}\right) =2u$ and $\delta \left( u\right) =0$. There is an
isomorphism $\phi :T\rightarrow S$ determined by $\phi \left( x_{1}\right)
=z_{1}+z_{2}$, $\phi \left( x_{2}\right) =z_{1}-z_{2}$, and $\phi \left(
t_{12}\right) $ is given by equation \ref{t form 2}. 
\begin{equation}
\phi \left( t_{12}\right) =\frac{1}{2}\left( 1-\varepsilon \left(
g_{1},g_{2}\right) \right) \left( \left( z_{1}\right) ^{2}-\left(
z_{2}\right) ^{2}\right) +\left( 1+\varepsilon \left( g_{1},g_{2}\right)
\right) \left( u-z_{1}z_{2}\right)  \label{t form 2}
\end{equation}%
\ Then $\phi \left( t\right) =u$ so $U\left( \mathcal{L}\right) \cong
T/\left( u\right) $. 
\endproof%

\begin{example}
\label{Color Example}Suppose $\varepsilon $ is a skew-symmetric bicharacter
on $G$ and there exist $g_{1},g_{2}\in G_{-}$ such that $2g_{1}=2g_{2}$, and 
$\varepsilon \left( g_{1},g_{2}\right) =1$\ (bicharacters with this property
are described in \cite[Lemma 2.7]{Prime}). Let $\mathcal{L}$ be the Lie
color algebra with homogeneous basis $\left\{
u_{1},u_{2},x_{1},x_{2}\right\} $ such that $\partial u_{1}=2g_{1}=2g_{2}$, $%
\partial u_{2}=g_{1}+g_{2}$, $\partial x_{1}=g_{1}$, $\partial x_{2}=g_{2}$,
and the brackets between basis elements are all zero except for the ones
listed below. 
\begin{equation*}
\begin{tabular}{cccc}
$\left\langle x_{1},x_{1}\right\rangle =2u_{1}$ & $\left\langle
x_{2},x_{2}\right\rangle =2u_{1}$ & $\left\langle x_{1},x_{2}\right\rangle
=u_{2}$ & $\left\langle x_{2},x_{1}\right\rangle =-u_{2}$%
\end{tabular}%
\end{equation*}%
Choose $\zeta \in k$ which satisfies $\zeta ^{2}=-1$, then $x_{1}+\zeta
x_{2} $ is a torsion element of $\mathcal{L}$. However $U\left( \mathcal{L}%
\right) \cong k\left[ \theta _{1}\right] \left[ \theta _{2};\sigma _{-1}%
\right] $ by the proof of Lemma \ref{Surprise}.
\end{example}

\begin{remark}
The product of linearly independent elements $x_{1}-x_{2},x_{1}+x_{2}\in 
\mathcal{L}_{-}$ is $u_{2}\in \mathcal{L}_{+}$. At first glance this may
appear to violate the PBW theorem (see \cite[Theorem 3.2.2]{BMPZ}). It does
not since the elements $x_{1}-x_{2}$ and $x_{1}+x_{2}$ are not homogeneous.
\end{remark}

\begin{remark}
In example \ref{Color Example}, $\mathcal{L}\cong L^{\gamma }$ for some
torsion free Lie superalgebra $L$ with appropriate $G$-grading and $\gamma $
a two-cocycle on $G$ (see \cite[Corollary 6.1]{Prime} for notation).\ There
is an algebra isomorphism $U\left( \mathcal{L}\right) \cong U\left( L\right)
^{\gamma }\cong U\left( L\right) $.
\end{remark}

\begin{theorem}
Let $\left( \mathcal{L},\left\langle ,\right\rangle \right) $ be a finite
dimensional Lie color algebra. Suppose $\dim \mathcal{L}_{-}\leq 2$ and $%
\dim V\leq \dim \left\langle V,V\right\rangle $ for each graded subspace $V$
of $\mathcal{L}_{-}$. Then $U\left( \mathcal{L}\right) $ is a domain with
global dimension equal to $\dim \mathcal{L}_{+}$.\label{Global Dim Thm}
\end{theorem}

\proof%
We may reduce to the case $\mathcal{L}_{+}$ is color central by \cite[Lemma
6.2]{Prime} and \cite[Corollary 6.18]{McRob}.\ Let $\mathcal{L}^{\prime }$
be the smallest sub Lie color algebra which contains $\mathcal{L}_{-}$. Then 
$\mathcal{L}_{+}^{\prime }=\left\langle \mathcal{L}_{-},\mathcal{L}%
_{-}\right\rangle $, $\mathcal{L}_{-}^{\prime }=\mathcal{L}_{-}$, and $%
U\left( \mathcal{L}\right) $ is a color polynomial algebra over $U\left( 
\mathcal{L}^{\prime }\right) $. Thus we may pass to the case $\mathcal{L}=%
\mathcal{L}^{\prime }$. If $\dim \mathcal{L}_{-}^{\prime }=\dim \mathcal{L}%
_{+}^{\prime }=2$ then we may apply Lemma \ref{Surprise}. Otherwise apply
Lemma \ref{Generic Ore}. 
\endproof%

\section{Gr\"{o}bner Basis Methods \label{Grob Method}}

Throughout this section $\mathcal{X}$ denotes a generic Lie color algebra
and $x_{1},x_{2},\ldots ,x_{m}$ form a homogeneous basis of $\mathcal{X}_{-}$%
. We define a chain of generic sub Lie color algebras which are used to
express $U\left( \mathcal{X}\right) $ as an iterated Ore extension. Then we
explain how to find the Gr\"{o}bner basis of an ideal generated by
homogeneous elements of $\mathcal{X}_{+}$.

Set $l=\frac{1}{2}m\left( m-1\right) $ and $p=l+m$. We specify a homogeneous
subset of $\mathcal{X}$ and define a function $\phi :\left\{ 1,2,\ldots
,p\right\} \rightarrow \left\{ 1,2\right\} $.\ Fix the following notation
for $1\leq i\leq j\leq m$.

\begin{itemize}
\item $s\left( i,j\right) =\frac{1}{2}j\left( j-1\right) +i$

\item $t_{s\left( j,j\right) }=x_{j}$

\item $\phi \left( s\left( j,j\right) \right) =2$

\item If $i<j$ then $t_{s\left( i,j\right) }=\left\langle
x_{i},x_{j}\right\rangle $ and $\phi \left( s\left( i,j\right) \right) =1$.
\end{itemize}

For each $\alpha =\left( \alpha _{1},\alpha _{2},\ldots ,\alpha _{p}\right)
\in \mathbb{N}^{p}$ we write $\mathbf{t}^{\alpha }=t_{1}^{\alpha
_{1}}t_{2}^{\alpha _{2}}\cdots t_{p}^{\alpha _{p}}\in U\left( \mathcal{X}%
\right) $. \ Set $e_{i}=\left( 0,\ldots ,\underset{\left( i\right) }{1}%
,\ldots ,0\right) \in \mathbb{N}^{p}$ for $1\leq i\leq p$. \ Then a
homogeneous basis for $\mathcal{X}_{+}$ is $\left\{ \left( t_{i}\right)
^{\phi \left( i\right) }:1\leq i\leq p\right\} =\left\{ \mathbf{t}^{\phi
\left( i\right) e_{i}}:1\leq i\leq p\right\} $.

Proceeding as in Remark \ref{generic subalgebra}, we let $\mathcal{X}_{i}$
be the generic sub Lie color algebra generated by $\left\{
x_{1},x_{2},\ldots ,x_{i}\right\} $.

\begin{lemma}
\label{Generic Ore 2}Set $U_{1}=k\left[ t_{1}\right] $, and for $j=2,\ldots
,m$, define $U_{j}$ recursively by equation \ref{recursive}. 
\begin{equation}
U_{j}=U_{j-1}\left[ t_{s\left( 1,j\right) };\sigma _{s\left( 1,j\right) }%
\right] \cdots \left[ t_{s\left( j-1,j\right) -1};\sigma _{s\left(
j-1,j\right) }\right] \left[ t_{s\left( j,j\right) };\sigma _{s\left(
j,j\right) },\delta _{s\left( j,j\right) }\right]  \label{recursive}
\end{equation}

\begin{enumerate}
\item Then $\mathcal{B}=\left\{ \mathbf{t}^{\alpha }:\alpha \in \mathbb{N}%
^{p}\right\} $ is a basis for $U\left( \mathcal{X}\right) $,

\item $U\left( \mathcal{X}_{i}\right) \cong U_{i}$ for $1\leq i\leq m$, and

\item there is a $\mathcal{B}$-admissible ordering $\preceq $ on $U\left( 
\mathcal{X}\right) $ such that $t_{1}\prec t_{2}\prec \cdots \prec t_{p}$.
\end{enumerate}
\end{lemma}

\proof%
Part 1 follows from the PBW\ Theorem. Part 2 can be proved in the same way
as \cite[Theorem 3.1]{Le Bruyn}. Part 3 follows from part 2 and \cite[%
Theorem 1.10]{Effective}.\ The maps $\sigma _{1},\sigma _{2},\ldots ,\sigma
_{p}$ and $\delta _{3},\delta _{6},\ldots ,\delta _{p}$ are determined by
the rules below.

\begin{description}
\item[(i)] $\sigma _{j}\left( t_{i}\right) =\varepsilon \left(
t_{j},t_{i}\right) t_{i}$ for $1\leq i<j\leq p$

\item[(ii)] $\delta _{s\left( j,j\right) }\left( t_{i}\right) =\left\langle
t_{s\left( j,j\right) },t_{i}\right\rangle $ for $1\leq j\leq m$ and $1\leq
i<s\left( j,j\right) $
\end{description}

\endproof%

We adopt notation and terminology from \cite{Effective}.

\begin{definition}
Choose nonzero $u\in U\left( \mathcal{X}\right) $. By part 1 of Lemma \ref%
{Generic Ore 2} the expression in equation \ref{stand rep} is unique with $%
c_{\alpha }\in k$ for each $\alpha \in \mathbb{N}^{p}$.%
\begin{equation}
u=\sum_{\alpha \in \mathbb{N}^{p}}c_{\alpha }\mathbf{t}^{\alpha }
\label{stand rep}
\end{equation}

\begin{enumerate}
\item Equation \ref{stand rep} is called the \emph{standard representation}
of $u$.

\item The \emph{Newton diagram} of $u$ is $\mathcal{N}\left( u\right)
=\left\{ \alpha \in \mathbb{N}^{p}:c_{\alpha }\neq 0\right\} $.

\item The \emph{exponent} of $u$ is $\exp \left( u\right) =\max_{\preceq }%
\mathcal{N}\left( u\right) $.

\item The \emph{leading coefficient} of $u$ is $lc\left( u\right) =c_{\exp
\left( u\right) }$.
\end{enumerate}
\end{definition}

In particular, if $u\in \mathcal{X}_{+}$ then $\mathcal{N}\left( u\right)
\subseteq \left\{ \phi \left( i\right) e_{i}:1\leq i\leq p\right\} $.

\begin{definition}
Let $u_{1},u_{2}\in U\left( \mathcal{X}\right) $ be given and set $\alpha
=\left( \alpha _{1},\alpha _{2},\ldots ,\alpha _{p}\right) =\exp \left(
u_{1}\right) $ and $\beta =\left( \beta _{1},\beta _{2},\ldots ,\beta
_{p}\right) =\exp \left( u_{2}\right) $. \ Let $\gamma =\left( \gamma
_{1},\gamma _{2},\ldots ,\gamma _{p}\right) $ be such that $\gamma _{i}=\max
\left\{ \alpha _{i},\beta _{i}\right\} $ for each $i=1,2,\ldots ,p$. \ The 
\emph{left }$S$\emph{-polynomial} of $u_{1}$ and $u_{2}$, denoted $S^{\ell
}\left( u_{1},u_{2}\right) $, is shown in equation \ref{S-polynomial} where $%
\lambda =lc\left( u_{2}\right) \left( lc\left( \mathbf{t}^{\alpha }\mathbf{t}%
^{\gamma -\alpha }\right) \right) ^{-1}$ and $\mu =lc\left( u_{1}\right)
\left( lc\left( \mathbf{t}^{\beta }\mathbf{t}^{\gamma -\beta }\right)
\right) ^{-1}$. 
\begin{equation}
S^{\ell }\left( u_{1},u_{2}\right) =\lambda \mathbf{t}^{\gamma -\alpha
}u_{1}-\mu \mathbf{t}^{\gamma -\beta }u_{2}  \label{S-polynomial}
\end{equation}
\end{definition}

\begin{lemma}
\label{Grobner Basis}If $K$ is a homogeneous linear subspace of $\mathcal{X}%
_{+}$ then there is a homogeneous basis $\mathcal{G}=\left\{
u_{1},u_{2},\ldots ,u_{n}\right\} $ of $K$ such that $u_{1}\prec u_{2}\prec
\cdots \prec u_{n}$ and $\exp \left( u_{i}\right) \notin \mathcal{N}\left(
u_{j}\right) $ for all $i,j$ with $1\leq i<j\leq n$.
\end{lemma}

\proof%
We prove such a basis $\mathcal{G}$ exists by induction on $n$, with the
case $n=1$ being trivial.

\begin{description}
\item[Step 1.] If $K^{\prime }$ is a homogeneous linear subspace of $K%
\mathcal{\ }$with $\dim K^{\prime }=n-1$ then there is a basis $\mathcal{G}%
^{\prime }=\left\{ u_{1},u_{2},\ldots ,u_{n-1}\right\} $ of $K^{\prime }$
such that $u_{1}\prec u_{2}\prec \cdots \prec u_{n-1}$ and $\exp \left(
u_{i}\right) \notin \mathcal{N}\left( u_{j}\right) $ for all $i,j$ with $%
1\leq i<j\leq n-1$.
\end{description}

This follows from the inductive hypothesis.

\begin{description}
\item[Step 2.] If $K^{\prime }$ and $G^{\prime }$ are as in Step 1 then
there exists homogeneous $u\in K\backslash K^{\prime }$ such that $\exp
\left( u\right) \neq \exp \left( u_{n-1}\right) $.
\end{description}

Suppose $v\in K\backslash K^{\prime }$ is homogeneous and $\exp \left(
v\right) =\exp \left( u_{n-1}\right) $. Then $\partial v=\partial u_{n-1}$
so $u=v-lc\left( v\right) \left( lc\left( u_{n-1}\right) \right)
^{-1}u_{n-1}\in K\backslash K^{\prime }$ is homogeneous with $\exp \left(
u\right) \prec \exp \left( u_{n-1}\right) $.

\begin{description}
\item[Step 3.] There exist $K^{\prime }$ and $\mathcal{G}^{\prime }$ as in
Step 1 such that some homogeneous $u\in K\backslash K^{\prime }$ satisfies $%
\exp \left( u_{n-1}\right) \prec \exp \left( u\right) $.
\end{description}

Let $K^{\prime }$ and $G^{\prime }$ be as in Step 1. By Step 2 there exists
homogeneous $w\in K\backslash K^{\prime }$ such that $\exp \left( w\right)
\neq \exp \left( u_{n-1}\right) $. If $\exp \left( w\right) \prec \exp
\left( u_{n-1}\right) $ set $u=u_{n-1}$ and let $K^{\prime \prime }$ be the
linear subspace spanned by $\left\{ u_{1},u_{2},\ldots ,u_{n-2},w\right\} $.
It is easy to see $\exp \left( w^{\prime }\right) \prec \exp \left(
u_{n-1}\right) $ for all $w^{\prime }\in K^{\prime \prime }$. Choose a basis 
$\mathcal{G}^{\prime \prime }$ as in Step 1 and replace $K^{\prime }$ with $%
K^{\prime \prime }$ and $\mathcal{G}^{\prime }$ with $\mathcal{G}^{\prime
\prime }$.

\begin{description}
\item[Step 4.] There is a homogeneous basis $\mathcal{G}=\left\{
u_{1},u_{2},\ldots ,u_{n}\right\} $ of $K$ such that $u_{1}\prec u_{2}\prec
\cdots \prec u_{n}$ and $\exp \left( u_{i}\right) \notin \mathcal{N}\left(
u_{j}\right) $ for all $i,j$ with $1\leq i<j\leq n$.
\end{description}

Let $K^{\prime }$, $\mathcal{G}^{\prime }$, and $u$ be as in Step 3 and
write $u$ as in equation \ref{stand rep}. Set $\mathcal{M}=\mathbb{N}%
^{p}\backslash \left\{ \exp \left( u_{i}\right) :i=1,2,\ldots ,n-1\right\} $
and set $u_{n}=u-\sum_{i=1}^{n-1}c_{\exp \left( u_{i}\right) }\left(
lc\left( u_{i}\right) \right) ^{-1}u_{i}$. Clearly $\exp \left(
u_{n-1}\right) \prec \exp \left( u_{n}\right) $ so we only need to show $%
\mathcal{N}\left( u_{n}\right) \subseteq \mathcal{M}$. Set $%
v=u-\sum_{i=1}^{n-1}c_{\exp \left( u_{i}\right) }\mathbf{t}^{\exp \left(
u_{i}\right) }$ and set $v_{i}=u_{i}-lc\left( u_{i}\right) \mathbf{t}^{\exp
\left( u_{i}\right) }$ for $i=1,2,\ldots ,n-1$. Then $\mathcal{N}\left(
v\right) \subseteq \mathcal{M}$ by construction and $\mathcal{N}\left(
v_{i}\right) \subseteq \mathcal{M}$ for $i=1,2,\ldots ,n-1$ by our choice of 
$\mathcal{G}^{\prime }$. Moreover $u_{n}=v-\sum_{i=1}^{n-1}c_{\exp \left(
u_{i}\right) }lc\left( u_{i}\right) ^{-1}v_{i}$ so $\mathcal{N}\left(
u_{n}\right) \subseteq \mathcal{N}\left( v\right) \cup \left(
\bigcup\limits_{i=1}^{n-1}\mathcal{N}\left( v_{i}\right) \right) \subseteq 
\mathcal{M}$ as desired.\ 
\endproof%

\begin{theorem}
\label{Grobner Ideal}Let $K$ be a homogeneous linear subspace of $\mathcal{X}%
_{+}$ and let $I$ be the ideal generated by $K$. Then $I=U\left( \mathcal{X}%
\right) K$ and $U\left( \mathcal{X}/K\right) \cong U\left( \mathcal{X}%
\right) /I$. Moreover the basis $\mathcal{G}$\ of $K$ described in Lemma \ref%
{Grobner Basis} is\ a Gr\"{o}bner basis of $I$.
\end{theorem}

\proof%
It is easy to see $\mathcal{X}_{+}$, and hence $K$, contains only normal
elements of $U\left( \mathcal{X}\right) $. The relations $I=U\left( \mathcal{%
X}\right) K$ and $U\left( \mathcal{X}/K\right) \cong U\left( \mathcal{X}%
\right) /I$ can be proved using the PBW\ Theorem (see \cite[Thoerem 3.2.2]%
{BMPZ}) and the universal property of the enveloping algebra.

To prove $\mathcal{G}$ is a Gr\"{o}bner basis of $I$ we show the remainder
is 0 when the left division algorithm by $\mathcal{G}$ is applied to $%
S^{\ell }\left( u_{i},u_{j}\right) $ (see \cite[Theorem 2.1]{Effective}).
Then $\mathcal{G}$ satisfies the so-called \textquotedblleft Buchberger
S-pair criterion,\textquotedblright\ which is \cite[Theorem 3.2]{Effective}.

For each $i,j$ with $1\leq i<j\leq n$ it is enough to show that $S^{\ell
}\left( u_{i},u_{j}\right) =v_{1}u_{i}+v_{2}u_{j}$ where $\exp \left(
v_{1}u_{i}\right) \preccurlyeq \exp \left( S^{\ell }\left(
u_{i},u_{j}\right) \right) $, $\exp \left( v_{2}u_{j}\right) \preccurlyeq
\exp \left( S^{\ell }\left( u_{i},u_{j}\right) \right) $ and for all $\alpha
\in \mathcal{N}\left( v_{i}\right) $ there does not exist $\beta \in \mathbb{%
N}^{p}$ such that $\exp \left( u_{j}\right) +\alpha =\exp \left(
u_{i}\right) +\beta $. We may pass to the case $i=1$ and $j=2$ by \cite[%
Proposition 2.11]{Effective}.

Since $u_{1},u_{2}\in \mathcal{X}_{+}$ and $u_{1}\prec u_{2}$ there exist $%
p_{1},p_{2}\in \left\{ 1,2,\ldots ,p\right\} $ with $p_{1}<p_{2}\leq p$ such
that $\exp \left( u_{1}\right) =\phi \left( p_{1}\right) e_{p_{1}}$ and $%
\exp \left( u_{2}\right) =\phi \left( p_{2}\right) e_{p_{2}}$. Thus $\alpha
=\phi \left( p_{1}\right) e_{p_{1}}$, $\beta =\phi \left( p_{2}\right)
e_{p_{2}}$, and $\gamma =\alpha +\beta $, which implies $\lambda =lc\left(
u_{2}\right) $ and $\mu =lc\left( u_{1}\right) \varepsilon \left(
u_{2},u_{1}\right) $. Set $v_{1}=\varepsilon \left( u_{2},u_{1}\right)
\left( u_{2}-lc\left( u_{2}\right) t_{p_{2}}^{\phi \left( p_{2}\right)
}\right) $ and $v_{2}=\varepsilon \left( u_{2},u_{1}\right) \left(
u_{1}-lc\left( u_{1}\right) t_{p_{1}}^{\phi \left( p_{1}\right) }\right) $.
Then%
\begin{eqnarray*}
S^{\ell }\left( u_{1},u_{2}\right) &=&\lambda \mathbf{t}^{\gamma -\alpha
}u_{1}-\mu \mathbf{t}^{\gamma -\beta }u_{2} \\
&=&lc\left( u_{2}\right) t_{p_{2}}^{\phi \left( p_{2}\right) }u_{1}-lc\left(
u_{1}\right) \varepsilon \left( u_{2},u_{1}\right) t_{p_{1}}^{\phi \left(
p_{1}\right) }u_{2} \\
&=&\left( u_{2}+v_{1}\right) u_{1}-\left( \varepsilon \left(
u_{2},u_{1}\right) u_{1}-v_{2}\right) u_{2} \\
&=&v_{1}u_{1}+v_{2}u_{2}
\end{eqnarray*}%
with $\exp \left( v_{1}\right) \prec \exp \left( u_{2}\right) $ and $\exp
\left( v_{2}\right) \prec \exp \left( u_{1}\right) $.

Let $\alpha \in \mathcal{N}\left( u_{2}\right) $ be arbitrarily chosen. Then 
$\alpha =\phi \left( p_{3}\right) e_{p_{3}}$ for some $p_{3}\in \left\{
1,2,\ldots ,p\right\} $ with $p_{3}\leq p_{2}$ since $u_{2}\in \mathcal{X}%
_{+}$. If there exists $\beta \in \mathbb{N}^{p}$ such that $\exp \left(
u_{j}\right) +\alpha =\exp \left( u_{i}\right) +\beta $ then $\phi \left(
p_{2}\right) e_{p_{2}}+\phi \left( p_{3}\right) e_{p_{3}}=\phi \left(
p_{1}\right) e_{p_{1}}+\beta $. We must have $\phi \left( p_{3}\right)
e_{p_{3}}=\phi \left( p_{1}\right) e_{p_{1}}$ and $\beta =\phi \left(
p_{2}\right) e_{p_{2}}$ since $p_{1}<p_{2}$. But this implies $\exp \left(
u_{1}\right) =\alpha \in \mathcal{N}\left( u_{2}\right) $, which contradicts
our choice of basis.

It follows from what we just proved that $\exp \left( v_{1}u_{1}\right) \neq
\exp \left( v_{2}u_{2}\right) $. This implies $\exp \left( v_{1}u_{1}\right)
\preccurlyeq \exp \left( S^{\ell }\left( u_{1},u_{2}\right) \right) $ and $%
\exp \left( v_{2}u_{2}\right) \preccurlyeq \exp \left( S^{\ell }\left(
u_{1},u_{2}\right) \right) $ as desired. 
\endproof%

\end{document}